\documentclass[12pt,]{amsart}
\usepackage{amsmath}
\usepackage{amscd}
\usepackage{amssymb}
\usepackage{latexsym}

\newtheorem{theorem}{Theorem}[section]
\newtheorem{lemma}[theorem]{Lemma}
\newtheorem{corollary}[theorem]{Corollary}
\newtheorem{proposition}[theorem]{Proposition}
\newtheorem{defprop}[theorem]{Proposition-Definition}
\newtheorem{remark}[theorem]{Remark}
\newtheorem{definition}[theorem]{Definition}

\newcommand{\bgl}{\begin{equation}}         
\newcommand{\egl}{\end{equation}}
\newcommand{\bgln}{\begin{eqnarray}}        
\newcommand{\egln}{\end{eqnarray}}
\newcommand{\bglnoz}{\begin{eqnarray*}}     
\newcommand{\eglnoz}{\end{eqnarray*}}
\newcommand{\btheo}{\begin{theorem}}
\newcommand{\etheo}{\end{theorem}}
\newcommand{\blemma}{\begin{lemma}}
\newcommand{\elemma}{\end{lemma}}
\newcommand{\bproof}{\begin{proof}}
\newcommand{\eproof}{\end{proof}}
\newcommand{\bbew}{\begin{beweis}}
\newcommand{\ebew}{\end{beweis}}
\newcommand{\bremark}{\begin{remark}\em}
\newcommand{\eremark}{\end{remark}}
\newcommand{\bdefin}{\begin{definition}}
\newcommand{\edefin}{\end{definition}}
\newcommand{\bprop}{\begin{proposition}}
\newcommand{\eprop}{\end{proposition}}
\newcommand{\bcor}{\begin{corollary}}
\newcommand{\ecor}{\end{corollary}}
\newcommand{\mn}{\par\medskip\noindent}

\newcommand{\cB}{\mathcal B}
\newcommand{\cC}{\mathcal C}

\newcommand{\cT}{\mathfrak T}

\newcommand{\lori}{\longrightarrow}
\newcommand{\lole}{\longleftarrow}

\newcommand{\ve}{\varepsilon}

\def\SEMI{\mbox{$\times\kern-2pt\vrule height5pt width.6pt \kern3pt $}}
\def\stackrellow#1#2{\mathrel{\mathop{#1}\limits_{#2}}}

\newcommand{\Hom}{{\rm Hom}\,}

\newcommand{\dis}{\displaystyle}

\def\Cz{\mathbb{C}}

\def\Zz{\mathbb{Z}}

\begin{document}
\title[Quillen's work on the foundations of cyclic cohomology]
{Quillen's work on the foundations of cyclic cohomology}
\date{\today}
\author[J. Cuntz]{Joachim Cuntz}

\address{Mathematisches Institut, Einsteinstr.62, 48149
M\"unster, Germany}\email{cuntz@uni-muenster.de}
\thanks{Research supported by DFG through
CRC 878 and by ERC through AdG 267079} 
\begin{abstract}\noindent
We give a survey of Quillen's contributions, apart from the very
first result in \cite{LQ}, to the area of cyclic homology.

\end{abstract}\maketitle

\begin{center}{Dedicated to the memory of Daniel
Quillen}\end{center}

\section{Introduction}

Daniel Quillen had been interested in cyclic homology from the very
start. After the first paper on the connection between cyclic
homology and the Lie algebra homology of matrices over an algebra
\cite{LQ}, his papers aimed at a deep understanding of the
fundamental structures underlying cyclic theory. When I went through
his articles on the subject I was again impressed by the systematic
development of his thoughts on the subject and by the very thorough
analysis of every detail which never left any loose ends.

A starting point for his work on the subject was the search for a
more conceptual explanation for Connes' construction of cyclic
cocycles from algebra extensions. He did this in \cite{CCAE},
\cite{CSFc} using an elegant formalism involving noncommutative
versions of Chern character forms and of Chern-Simon forms. He also
showed that the important JLO-cocycle \cite{JLO} could be obtained
from a similar construction. In \cite{CCAE}, he showed that all
cocycles for an algebra $A$ can be obtained from an extension of the
form $I\rightarrowtail R\twoheadrightarrow A$ using traces on
$R/I^n$ or on $I^n$. In this way he could represent the cyclic
cohomology groups $HC^{2n}A$ and $HC^{2n+1}A$, as inductive limits,
over all extensions of $A$, of the spaces of traces on $R/I^{n+1}$
or on $I^{n+1}$, respectively.

At about the same time, A.Connes and the author had studied traces
on the ideal $\ve A$ in a universal so-called semi-split extension
$\ve A\rightarrowtail EA\twoheadrightarrow A$ and had shown that one
obtains the cyclic cohomology from traces on powers of that ideal
(or quotients by these powers) \cite{CoCu}. It was noted in
\cite{CuUE} that there is a bijection between traces on powers of
$\ve A$ and traces on powers of the ideal in a free extension of
$A$. This observation gave another proof of Quillen's result in
\cite{CCAE} and was the beginning of a long lasting collaboration of
Dan Quillen and the author. This collaboration culminated in
\cite{CQCNs} and in a new approach to cyclic homology (especially
well adapted to the periodic theory) and in a proof of excision in
periodic cyclic theory in \cite{CQEx}. Some of the concepts
developed along the way such as in particular the one of smoothness
(quasi-freeness) in noncommutative geometry also became influential.

A basis for the approach to cyclic homology developed in
\cite{CQCNs} is a very simple - nearly trivial - periodic complex of
period 2, the $X$-complex $X(A)$ which was already used by Quillen
in \cite{ACCC}. It is a reduction of the cyclic bicomplex (or
equivalently of the $B,b$-bicomplex) to the lowest dimensions 0 and
1. Even though it is so simple it is still good enough to capture
the cyclic (and Hochschild) homology for algebras of homological
dimension $\leq 1$. Since, in particular, free algebras have that
property, such algebras were called quasi-free in \cite{CQANs}. The
properties of quasifree algebras were analyzed in \cite{CQANs}.
Interestingly they turn out to be exactly the natural generalization
of the notion of a smooth variety or algebra, to the non-commutative
setting. In fact, they are characterized by a lifting property which
is the exact analogue of the corresponding property of a smooth
algebra in the commutative category.

For an arbitrary algebra $A$ the cyclic homology/cohomology can then
be obtained in the following way. Choose any extension of $A$ of the
form $0\to I\to R\to A\to 0$ with $R$ quasifree and consider the
$I$-adic completion $\hat{R}=\mathop{\lim}\limits_{
{\scriptstyle\longleftarrow} } R/I^n$.

Then the periodic cyclic homology $HP_*(A)$ is simply the homology
of the complex $X_*(\hat{R})$, the periodic cyclic cohomology
$HP^*(A)$ is the homology of the (continuous for the $I$-adic
topology) dual of $X_*(A)$ and the ordinary cyclic
homology/cohomology groups $HC_n(A), HC^n(A)$ can be obtained from a
natural filtration (a small modification of the $I$-adic filtration)
of these complexes. This procedure is directly analogous to the
construction of infinitesimal homology in algebraic geometry where
one embeds a general variety into a smooth variety, completes and
the considers the de Rham complex of the completion. The $X$-complex
thus plays the role of a noncommutative de Rham complex.

The simplicity of this description of cyclic theory is a little bit
obscured in \cite{CQCNs} by the fact that the exposition there
strives to elucidate all facets of the approach as well as its
connections to other approaches.

The concepts involved in this approach also were a natural basis for
an attack on the problem of excision for periodic cyclic
homology/cohomology. The excision problem for the ordinary cyclic
theory had been understood for quite a while. Wodzicki \cite{Wod}
had shown that excision does not always hold and that it holds,
given an algebra $I$, for any extension of the form
$I\rightarrowtail A\twoheadrightarrow A/I$ if and only if $I$ has a
property which he called $H$-unitality. On the other hand,
Goodwillie \cite{Good} had shown that excision in the periodic
cyclic theory also holds for nilpotent ideals (which are never
$H$-unital). It was therefore natural to ask whether excision holds
for arbitrary extensions in the periodic cyclic theory. That problem
had remained open for a long time. In \cite{CQCR2} it was discovered
that every ideal in a quasi-free algebra satisfies a property which
was called approximate $H$-unitality and that this property is
enough to prove excision. Since any algebra can be represented as a
quotient $R/I$ of a quasi-free (even free) algebra and since the
periodic cyclic cohomology of $A$ can be easily related to that of
$I$ in such an extension, excision in periodic cyclic cohomology
followed in complete generality.

The proof of excision for periodic homology (rather than cohomology)
and for the bivariant case needed some additional ideas and
techniques. It fits naturally into the framework of the $X$-complex
and quasi-free extensions, and was developed in \cite{CQEx}. The
proof in that paper takes a detour by reducing the general problem
first to the quasi-free case and then applying an argument \`{a} la
Wodzicki. This detour was avoided in a later simplified proof due to
R.Meyer \cite{Meyexc}, which fits perfectly into the framework of
\cite{CQCNs}.

In this note we give a brief introduction to some of the main ideas in the series of 10 articles that have been written in the period between 1987 and 1997. The reader will notice the linear progression of Quillen's thoughts on the subject. Several of the results from that period represent fundamental new findings. The notion of noncommutative smoothness (being quasi-free) that has been analyzed in much detail in \cite{CQANs}, has become influential in various contexts, see e.g. \cite{KR}. The new description of cyclic theory in \cite{CQCNs} has become the framework of choice for topological theories such as entire and local cyclic theory, cf. \cite{Meybook}, \cite{PuLoc}, but also for equivariant cyclic theory, \cite{Voi1}, \cite{Voi2}.
Finally, the excision result \cite{CQEx} has opened the way to a better understanding of cyclic theory, to the construction of a bivariant Chern character on suitable categories of algebras,  \cite{CuEnz} \cite{PuLoc} and has been generalized to other cyclic theories such as the entire and the local theory.

\section{Cyclic homology and algebra extensions}

Let $I\rightarrowtail R\twoheadrightarrow A$ be an extension (i.e. a short exact sequence
where the arrows are algebra homomorphisms) of the algebra $A$. Note that, in this article, every algebra is an algebra over a field of characteristic 0, which we usually even assume to be the field of complex numbers $\Cz$. To a
trace (i.e. a linear functional vanishing on commutators) on $R/I^n$
or on $I^{n+1}$, Connes \cite{CoNCG} had associated cyclic cocycles
in $HC^{2n}$ and $HC^{2n+1}$, respectively. Dually, Connes'
construction leads to natural maps

\bgl\label{map}HC_{2n}(A)\to\; R\,\big/(I^{n+1}+[R,R])\qquad
HC_{2n+1}(A)\to\; I^{n+1}\big/[I,I^n]\egl

where $[\cdot\,,\,\cdot]$ denotes the linear space generated by all
commutators. In \cite{ACCC}, Quillen, introducing an elegant
formalism, showed how to interpret the cyclic cocycles obtained
using the dual maps from the spaces of traces on $R/I^n$ or on
$I^{n+1}$ to cyclic cohomology, as Chern character maps and
Chern-Simons forms. He also showed that for large $n$ they are
related by Connes' $S$-operator. For these considerations he used
already formulas that came to be important later in \cite{CQCNs} and
in particular he already used the $X$-complex in order to prove the
$S$-relations. His construction of cyclic cocycles using the
formalism of Chern and Chern-Simons forms has been used later by
various authors in connection with index theorems, see e.g.
\cite{Per}, \cite{Hig}.\mn

In \cite{CCAE}, Quillen then showed that the maps in (\ref{map}) are
injective, if $R$ is a free algebra. Thus, every cyclic cocycle for
$A$ can be represented by a trace on $R/I^{n}$ or on $I^n$, if $R$
is a free extension of $A$. From this he then obtains the following
description of cyclic homology

$$HC_{2n} = \mathop{\lim}\limits_{
{\scriptstyle\longleftarrow} }
\left(R\Big/(I^{n+1}+[R,R])\right)\qquad HC_{2n+1} =
\mathop{\lim}\limits_{ {\scriptstyle\longleftarrow} }
\left(I^{n+1}\Big/[I,I^n]\right)$$

where the inverse limit is taken over all extensions of $A$ of the
form $0 \to I \to R \to A \to 0$.

On the other hand, at about the same time when Quillen was working
on \cite{CCAE}, Connes and the author had analyzed traces on the
free product algebra $QA =A\star A$ and had shown that they are
described by families of multilinear functionals that describe
cocycles in the cyclic bicomplex \cite{CoCu}. It was noted in
\cite{CuUE} that the canonical free extension $RA$ of an algebra $A$
given by the tensor algebra over $A$ is a natural subalgebra of $QA$ (the
even part for a natural $\Zz/2$-grading) and that therefore, traces
on $RA$ correspond bijectively to graded traces on $QA$. This
observation gave an alternative proof for Quillen's results in
\cite{CCAE} and was the starting point for a long lasting cooperation
between the author and Dan Quillen.

At the end of the introduction to the paper \cite{CCAE}, Quillen
writes: ``It is clear from the present paper, with its extensive use
of explicit complexes and formulas, that a true Grothendieck
understanding of cyclic cohomology remains a goal for the future.
Indeed, the inverse limit formula for cyclic homology described
above is a strong indication that there is a much simpler foundation
of the subject.''

The search for such a formulation was the guideline in Quillen's
subsequent work and in our collaboration.

\section{Operators on differential forms and convenient description
of the cyclic bicomplex}\label{P}

We will have to refer to one of the standard definitions of cyclic
homology. The most convenient way to describe the fundamental cyclic
bicomplex for our purposes is in the guise of the $B,b$-bicomplex
defined by operators $B$ and $b$ on the algebra $\Omega A$ of
differential forms over a given algebra $A$. We go immediately
\textit{in medias res} and introduce the important harmonic
decomposition of $\Omega A$ which has been developed in \cite{CQOp}
and used in \cite{CQCNs} by Cuntz-Quillen.

Given an algebra $A$, we denote by $\Omega A$ the universal algebra
generated by $x\in A$ with relations of $A$ and symbols $dx,x\in A$,
where $dx$ is linear in $x$ and satisfies $d(xy)=xd(y) + d(x)y$. We
do not impose $d1=0$, i.e., if $A$ has a unit, $d1\neq0$. $\Omega A$
is a direct sum of subspaces $\Omega^nA$ generated by linear
combinations of $\; x_0dx_1 \,\dots\, dx_n\,$, and $\; dx_1
\,\dots\, dx_n\, ,\; x_j\in A$. This decomposition makes $\Omega A$
into a graded algebra. We write $\deg(\omega)=n$ if
$\omega \in \Omega^nA$. \mn As a vector space, for $n\ge 1$,

\bgl\label{iso1} \Omega^nA \cong \widetilde A \otimes A^{\otimes
n}\cong A^{\otimes (n+1)} \oplus A^{\otimes n} \egl

(where $\widetilde A$ is $A$ with a unit adjoined, and $1\otimes
x_1\otimes\dots\otimes x_n$ corresponds to $ dx_1\dots dx_n$). The
operator $d$ is defined on $\Omega A$ by
$$
d(x_0dx_1\dots dx_n)\, =\, dx_0dx_1\dots dx_n \qquad d(dx_1\dots
dx_n)\, =\, 0
$$
The operator $b$ is defined by
$$
b(\omega dx) = (-1)^{\deg\omega}[\omega,x] \qquad b(dx) = 0\, ,\,
b(x) \, = \, 0\, , \qquad x\in A\, ,\, \omega\in\Omega A
$$
Then clearly $d^2 = 0$ and one easily computes that also $b^2 = 0$.
Under the isomorphism in equation (\ref{iso1}) $d$ becomes
$$
d(x_0 \otimes \ldots \otimes x_n) = \mbox{1}\otimes x_0 \otimes
\ldots \otimes x_n\qquad d(\mbox{1} \otimes x_1 \otimes \ldots
\otimes x_n) = 0
$$
while $b$ corresponds to the usual Hochschild operator $A^{\otimes
(n+1)}\to A^{\otimes n}$. Another important natural operator is the
degree (or number) operator defined by $ N(\omega)\;  =
\;\;\mbox{deg}(\omega) \omega$.
Now, for the operator  $L=(Nd)b +b(Nd)$ one obtains a splitting
$\Omega A=\mbox{Ker}\,L\oplus\mbox{Im}\,L$ (this follows from a polynomial identity satisfied by the
Karoubi operator $\kappa=db+bd$).

The operator $L$ thus behaves like a ``selfadjoint" operator. It can
be viewed as an abstract Laplace operator on the algebra of abstract
differential forms $\Omega A$. We denote by $P$ the projection onto the kernel of $L$. The elements
in the image of $P$ are then ``abstract harmonic forms", \cite{CQOp}, \cite{CQCNs}. \mn By
construction, $P$ commutes with $b$, $d$, $N$. Thus setting $B=NPd$
one finds $Bb+bB=PL=0$ and $B^2=0$.

The preceding identities show that we obtain a bicomplex - the
$(B,b)$-bicomplex - in the following way

\bgl\label{Bb}
\begin{array}{ccccccc}
\downarrow b & & \downarrow b & & \downarrow b & &\downarrow b \\[0.2cm]
\Omega^3A & \stackrel{B}{\longleftarrow} & \Omega^2A &
\stackrel{B}{\longleftarrow} &
\Omega^1A & \stackrel{B}{\longleftarrow} & \Omega^0A\\[0.2cm]
\downarrow b & & \downarrow b & & \downarrow b & & \\[0.2cm]
\Omega^2A & \stackrel{B}{\longleftarrow} & \Omega^1A &
\stackrel{B}{\longleftarrow} &
\Omega^0 A &  & \\[0.2cm]
\downarrow b & & \downarrow b & &  & & \\[0.2cm]
\Omega^1A & \stackrel{B}{\longleftarrow} & \Omega^0 A &  & &  & \\[0.2cm]
\downarrow b & & & &  & & \\[0.2cm]
\Omega^0 A &&&&&&
\end{array}
\egl One can rewrite the $(B,b)$-bicomplex (\ref{Bb}) using the
isomorphism $\Omega^nA \cong A^{\otimes (n+1)} \oplus A^{\otimes n}$
in equation (\ref{iso1}). An easy computation shows that then it
becomes the usual cyclic bicomplex with the operators $b,b'$ and
Connes' signed cyclic permutation operator $\lambda$.

\bdefin One defines the cyclic homology $HC_n(A)$ of the algebra $A$
to be the homology of the total complex of the
$B,b$-bicomplex.\edefin

We denote by $\widehat{\Omega}A$ the infinite product
$$
\widehat{\Omega}A = \prod_n \Omega^n A
$$
and by $\widehat{\Omega}^{ev}A$, $\widehat{\Omega}^{odd}A$ its even
and odd part, respectively. $\widehat{\Omega}A$ may be viewed as the
(periodic) total complex for the bicomplex (\ref{Bb}) continued
infinitely to the left and down. Similarly, the (continuous for the
filtration topology) dual $(\widehat{\Omega}A)'$ of
$\widehat{\Omega}A$ is
$$
(\widehat{\Omega}A)' = \bigoplus_n \,(\Omega^n A)'
$$
\begin{definition}
The periodic cyclic homology $HP_*(A)$, $*=0,1$, is defined as the
homology of the $\Zz/2$-graded complex
\[ \widehat{\Omega}^{ev}A \quad
\stackrel{\dis\stackrel{B+b}{\longrightarrow}}
{\dis\stackrellow{\longleftarrow}{B+b}} \quad
\widehat{\Omega}^{odd}A
\]
and the periodic cyclic cohomology $HP^*(A)$, $*=0,1$, is defined as
the homology of the $\Zz/2$-graded complex
\[ (\widehat{\Omega}^{ev}A)' \quad
\stackrel{\dis\stackrel{B+b}{\longleftarrow}}
{\dis\stackrellow{\longrightarrow}{B+b}} \quad
(\widehat{\Omega}^{odd}A)'
\]
\end{definition}

The comparison of the $B,b$-bicomplex with the $X$-complex
for a quasi-free resolution suggests in fact the boundary operator
$B-b$ rather than $B+b$.. This convention would also avoid
complicated signs in the Chern character map. We follow here however
the conventions in \cite{CQCNs}.

\section{The $X$-complex and quasi-free algebras}

The $X$-complex is the quotient of the $B,b$-bicomplex (\ref{Bb}) by
the following sub-bicomplex

\bgl\label{sub}
\begin{array}{cccccc}
&\downarrow b & & \downarrow b & & \downarrow b \\[0.2cm]
\lole&\Omega^3A & \stackrel{B}{\longleftarrow} & \Omega^2A &
\stackrel{B}{\longleftarrow} &
b(\Omega^2 A) \\[0.2cm]
&\downarrow b & & \downarrow b & &  \\[0.2cm]
\lole&\Omega^2A & \stackrel{B}{\longleftarrow} & b(\Omega^2 A) & & \\[0.2cm]
&\downarrow b & & & &  \\[0.2cm]
\lole&b(\Omega^2A) &&&&
\end{array}
\egl

Thus $X(A)$ is the periodic complex

$$ X(A):\quad\to A \stackrel{\natural d}{\longrightarrow}
\Omega^1A_\natural \stackrel{b}{\longrightarrow} A
\stackrel{\natural d}{\longrightarrow} \Omega^1A_\natural
\to\quad$$

where $\Omega^1A_{\natural} = \Omega^1A/[A,\Omega^1A]$ is the
quotient of the bimodule $\Omega^1A$ by the subspace of commutators
and $\natural : \Omega^1A \to \Omega^1 A_{\natural}$ is the
canonical quotient map.

Even though the $X$-complex is very simple, it is good enough to
compute the periodic cyclic theory for a special class of
algebras - the quasi-free algebras for which
the subcomplex by which we divide is contractible.
\begin{defprop}\label{qfree} (\cite{Shel}, \cite{CQANs}) Let $A$ be an
algebra. The following conditions are equivalent:
\begin{enumerate}
\item[(1)] Let $0 \rightarrow N \rightarrow  S \stackrel{q}{\rightarrow}
B \rightarrow 0$ be an extension of algebras where the ideal $N$ is
nilpotent (i.e., $N^k = \{0\}$ for some $k \ge 1$) and $A
\stackrel{\alpha}{\longrightarrow} B$ a homomorphism. Then there
exists a homomorphism $A \stackrel{\alpha'}{\longrightarrow} S$ such
that $q \circ \alpha' = \alpha$.
\item[(2)] $A$ has cohomological dimension $\le 1$ with respect to Hochschild
cohomology.
\end{enumerate} The algebra  $A$ is called
quasi-free if these equivalent conditions are satisfied.
\end{defprop}

Many other different characterizations of quasi-freeness are analyzed in \cite{CQANs}. An especially important property for us is the fact that $HX_*(A) = HP_*A$ for a quasi-free algebra $A$ (where $HX_*(A), * = 0,1$ denotes the homology of the complex $X(A)$).\mn

The fact that $HP_*$ is invariant under polynomial (or differentiable) homotopies, is reflected by the fact that the $X$-complex has, for quasifree algebras a homotopy invariance
property which is described by a natural Cartan homotopy formula.

The most important examples of quasi-free algebras are free
algebras. In particular the non-unital tensor algebra $TA$ given by

\bgl\label{tens} TA = A\oplus A^{\otimes^2}\oplus
A^{\otimes^3}\oplus\ldots\egl

is quasi-free and the natural quotient map $TA\to A$ defines a
quasi-free extension of $A$.

Another interesting feature of quasi-free algebras, which justifies to consider quasi-freeness as the correct analogue of smoothness in the noncommutative situation, is the following ``tubular neighbourhood theorem'' proved in \cite{CQANs}.

\btheo\cite{CQANs}[6, Theorem 2]. Assume $A$ is quasi-free. Consider
an extension $A=R/I$ and let $N$ denote the $A$-bimodule $I/I^2$. If
$A$ and $R$ are quasi-free, there is an isomorphism
$u:\hat{T}_AN\to\hat{R}$ from the $N$-adic completion $\hat{T}_AN$
of the tensor algebra $T_AN$ to the $I$-adic completion of $R$ which
extends the identity map on $A\oplus N$.\etheo

The concept of smoothness for a noncommutative algebra studied in \cite{CQANs}, was later found to be great significance in noncommutative geometry, cf. e.g. \cite{KR}.

\section{Cyclic homology and nonsingularity}

Let $p : A \to TA$ be the canonical linear inclusion of $A$ into the
tensor algebra over $A$ (see (\ref{tens}) above) and, for $x,y \in
A$ set $\omega(x,y) = p(xy) - pxpy$. Then the map
\[
\alpha : x_0 dx_1 dx_2 \ldots dx_{2n-1} dx_{x_{2n}} \; \mapsto \;
x_0 \omega (x_1,x_2) \ldots \omega (x_{2n-1}, x_{2n})
\] defines a linear isomorphism $\Omega^{ev}A
\stackrel{\alpha}{\longrightarrow} TA$. In fact, this isomorphism
becomes an algebra isomorphism for a deformation of the product on
$\Omega^{ev}A$, the Fedosov product. It extends to a linear
isomorphism $\alpha: \Omega A\to X(TA)$ respecting the
$\Zz/2$-grading.

May be the most important theorem in \cite{CQCNs} which nearly seems
like a miracle is the following.

\btheo\label{isom} \cite{CQCNs} The map $\alpha$ gives an isomorphism (not just a quasi-isomorphism!) between the
$\Zz/2$-graded complexes $P\Omega A$ (with boundary operator $B+b$)
and $PX(TA)$ which respects a natural filtration on both
sides.\etheo

Here $P$ denotes the harmonic projection operator, see section
\ref{P}, on $\Omega A$ and also its counterpart on $X(TA)$. Due to
the definition and properties of $P$, the complements $P^\perp\Omega
A$ and $P^\perp X(TA)$ do not have any homology. As a consequence
all cyclic invariants such as cyclic homology, cyclic cohomology,
periodic cyclic theory, bivariant cyclic theory can be defined and
computed alternatively from the $B,b$-bicomplex or from the complex
$X(TA)$.

Thus one can define all cyclic
homology/cohomology-invariants using the simple complex $X(TA)$. In
fact, by homotopy invariance of the $X$-complex for quasi-free
algebras, one may even replace here $TA$ by any quasi-free extension
$T$ of $A$.

For instance we get for the periodic cyclic homology of $A$ the
formula

$$HP_*(A) = HX_*(\hat{T})$$

where $J\rightarrowtail T\twoheadrightarrow A$ is any extension of $A$ with $T$
quasi-free and $\hat{T} =\mathop{\lim}\limits_{
{\scriptstyle\longleftarrow} } T/J^n$ is the $J$-adic completion.
The ordinary cyclic homology groups $HC_nA$ can be determined from a
natural filtration of $X(T)$ (a small modification of the $J$-adic
filtration giving the same completion). The cohomology groups and
more generally, the bivariant periodic cyclic homology groups for
two algebras $A$ and $B$, are determined by the formula

$$ HP_*(A,B)= H_*(\Hom (X(\hat{T}),X(\hat{S}))$$

where $S$ is a completed quasi-free extension of $B$ and Hom denotes
the Hom-complex based on maps which are continuous for the adic
topologies. From this definition we immediately get a product $HP_i(A,B)\times HP_j(B,C)\to HP_{i+j}(A,C)$. For the periodic cyclic theory this definition of the bivariant theory, given in \cite{CQCNs} is better suited than the former definition of Jones-Kassel \cite{JK}. For instance any surjective algebra homomorphism $A\to B$ with nilpotent kernel induces an \emph{invertible} element in $HP_*(A,B)$. However the inverse can not be realized as an element of the Jones-Kassel bivariant theory $HC^*(B,A)$. Also, excision holds for $HP_*$ (see section \ref{exc}), but not for the Jones-Kassel theory.

One virtue of this new description of cyclic theory is that it
reduces many computations to computations with the very simple
$X$-complex for quasi-free algebras.

It has been shown in \cite{CQANs} that the fundamental properties of
homotopy invariance and invariance \cite{Good} under nilpotent
extensions, of the periodic theory, follow very naturally from this
description. In fact, the invariance under nilpotent extensions is
built in into the definition and homotopy invariance follows from a
Cartan homotopy formula for the $X$-complex.

\section{Models for cyclic homology types}

Motivated by the new description of cyclic homology in \cite{CQCNs},
Quillen in \cite{BCCM} went on to systematically compare the different ways that
had been used in the literature to describe cyclic homology.

A mixed complex is a complex with differential $B$ and equipped with
a second operator $B$ of degree +1 such that $[b,B]=B^2=0$. It is
called free when its homology with respect to the differential $B$
is zero.

An $S$-module is a complex with an operator $S$ of degree -2
commuting with the differential.  It is called divisible when the
operator $S$ is surjective.

A supercomplex is a $\Zz/2$-graded complex, i.e. a $\Zz/2$-graded
vector space with an odd operator $d$ of square 0. A tower of
supercomplexes is an inverse system $X=X^n$ of supercomplexes
indexed by the integers such that the maps $X^n\to X^{n-1}$ are
surjective supercomplex maps. It is called special if the odd degree
homology of the associated graded vanishes.

All mixed complexes, $S$-modules and towers of supercomplexes are
assumed to be bounded below in the sense that the components $C^n$
or $X^n$ are zero for $n\ll 0$.

Let $\cC_\Lambda$, $\cC_S^d$ and $\cT^s$ denote the categories of mixed
complexes, divisible $S$-modules and special towers, respectively.

There are natural comparison functors
$\cC_\Lambda\stackrel{\cB}{\lori}\cC_S^d\stackrel{\alpha}{\lori}\cT^s$.
Moreover there are Hochschild, cyclic and periodic cyclic homology functors
from each of these categories to vector spaces, and $\cB$, $\alpha$
are compatible with those functors.

One of the essential results in \cite{CQCNs} can be stated as saying
that the special tower obtained from a natural filtration of $X(R)$
for a quasi-free extension $I\rightarrowtail R\twoheadrightarrow A$ of the algebra A, is
homotopy equivalent to the image under $\alpha\circ\cB$ of the mixed
complex $\Omega A$ with operators $b$ and $B$.

Quillen shows that the following five categories are equivalent in
general:
\begin{enumerate}
  \item The derived category of mixed complexes.
  \item The homotopy category of free mixed complexes.
  \item The derived category of $S$-modules.
  \item The homotopy category of divisible $S$-modules.
  \item The homotopy category of special towers of supercomplexes.
\end{enumerate}
The morphisms in each category can be viewed as cycles in a
bivariant cyclic cohomology similar to the one of Jones-Kassel
\cite{JK}. Each object in any of the categories describes a cyclic
homology type in the sense that it contains the full information on
the cyclic homology, Hochschild homology, negative cyclic homology
and periodic cyclic homology. By the universal coefficient theorem
established by Jones-Kassel for their bivariant theory, a cyclic
homology type is determined, inside each of these categories, by its
cyclic homology groups.

As a byproduct of these considerations, Quillen also relates the
results in \cite{CQCNs} to his earlier results in \cite{CCAE}.

\section{Excision in periodic cyclic homology}\label{exc}
The results in \cite{CQCNs} had established a very satisfactory
framework for cyclic theory explaining smoothly some fundamental
properties such as homotopy invariance, Morita invariance,
invariance under nilpotent extensions and furnishing a natural
picture for the Chern character maps from $K$-theory to cyclic
homology. Especially in the periodic case the bivariant theory thus
now looked formally exactly like the equally $\Zz/2$-graded
bivariant topological $K$-theory of Kasparov which plays a fundamental role in
noncommutative geometry. It became completely clear that the missing
ingredient to relate $K$-theory and cyclic theory was excision - the
fact that any extension of algebras induces long exact sequences in
the two variables of the bifunctor.

The starting point for our attack on the excision problem in
periodic cyclic theory was the observation in \cite{CQCR1}, \cite{CQCR2} that if $I$ is an ideal
in a quasi-free algebra $R$, then the projective system $(I^n)$
of powers of $I$ satisfies a property which is analogous to Wodzicki's
$H$-unitality (namely the property that the multiplication map $I\otimes I\to I^2$ admits an $I$-linear splitting). Let $C(A)$ be the complex describing periodic cyclic
cohomology $HP^*A$ of an algebra $A$ and consider an extension of algebras
of the form $0\to I\to P\to Q\to 0$. We may then take the inductive
limit over $n$ for the extensions of complexes

$$0\to C(I^n:P)\to C(P)\stackrel{\pi}{\lori} C(Q/I^n)\to 0$$

where $C(I^n:P)$ denotes the relative complex, i.e. the kernel of
$\pi$. An argument \`{a} la Wodzicki shows that the cohomology of
the first complex will converge in the inductive limit to $HP^*(I)$
while, by invariance of $HP^*$ under nilpotent extensions
\cite{Good}, the cohomology of the third complex will converge to
$HP^*(Q)$. The limit of the associated long exact cohomology
sequence for this exact sequence of complexes will then be of the
form

$$
\begin{array}{cccccr}
  HP^0(I) & \leftarrow & HP^0(P) & \leftarrow & HP^0(Q)
  \vspace{1mm}& \\
  \downarrow & & & & \uparrow\vspace{1mm}&\qquad\\
  HP^1(Q) & \to & HP^1(P) & \to & HP^1(I)&
\end{array}
  $$

In a second step it was noted in \cite{CQCR2}, that every extension
$0\to I\to A\to B\to 0$ can be mimicked by an equivalent extension
$0\to I'\to A'\to B'\to 0$ where $I'$ is isomorphic to an ideal in a
quasi-free algebra. This established excision for periodic cyclic
\emph{cohomology} in  general. The excision problem for periodic
cyclic \emph{homology} remained open at that point.

It is possible, even though not trivial, to derive from
the excision result for periodic cyclic cohomology the one for
periodic cyclic homology, using homological algebra. A problem in
doing so is the more delicate behaviour of homology and exact
sequences under inverse limits rather than direct limits. However, such a
computational approach is not good enough. For one thing, it would not work for
topological algebras, but, more importantly, it would not give the
result for the bivariant theory. The excision result for the
bivariant theory is very important, even if one is only interested
in the monovariant theory. For instance, the connecting maps in the
long exact sequences for periodic cyclic homology and cohomology are
given by the product by a natural bivariant element determined by
the extension.

The method to attack that technical problem, that was used in
\cite{CQEx}, was to cast the entire theory into the framework of
pro-vector spaces, pro-complexes and at the same time to extend it
from algebras to pro-algebras (by definition a pro-object is an
inverse system in a category). An important feature of the category
of pro-vector spaces is the fact that it contains enough projective
objects. On the other hand, if $0\to I\to R\to A\to 0$ is a
quasi-free extension, then the pro-complex $X(R/I^n)$ which
describes the periodic theory of $A$ can be replaced by the
pro-complex $X(I^n)$ which  describes the same theory (up to a dimension
shift even $\leftrightarrow$ odd). However this second complex is a
complex of projective pro-vector spaces. Thus
Hom$(X(I^n),\,\,\cdot\,)$ is exact on pro-vector spaces.

Using this technology, the argument for excision in periodic cyclic
cohomology outlined in \cite{CQCR1} carries over to the bivariant theory - but in a
more conceptual way, since inductive and projective limits are
already incorporated into the notion of a morphism between
pro-vector spaces.

\begin{theorem}\label{excision} \cite{CQEx} Let $0 \rightarrow S \rightarrow P \rightarrow Q \rightarrow 0$ be an
extension of algebras and $A$ an algebra. There are two natural
six-term exact sequences
  $$
\begin{array}{cccccr}
  HP_0(A,S) & \to & HP_0(A,P) & \to & HP_0(A,Q)
  \vspace{1mm}& \\
  \uparrow & & & & \downarrow\vspace{1mm}&\qquad\\
  HP_1(A,Q) & \leftarrow & HP_1(A,P) & \leftarrow & HP_1(A,S)&
\end{array}
  $$
  and
  $$
\begin{array}{cccccr}
  HP_0(S,A) & \leftarrow & HP_0(P,A) & \leftarrow & HP_0(Q,A)
  \vspace{1mm}& \\
  \downarrow & & & & \uparrow\vspace{1mm}&\qquad\\
  HP_1(Q,A) & \to & HP_1(P,A) & \to & HP_1(S,A)&
\end{array}
  $$  where the horizontal arrows are induced by the maps in the given extension and the vertical maps are given by the product with a natural class in $HP_1(Q,S)$, associated with the extension.
\end{theorem}

A simplified proof of this theorem - even more directly in the
spirit of the framework of \cite{CQCNs} and avoiding Wodzicki's
argument completely - was given later by R.Meyer \cite{Meyexc}. Let
$S \rightarrowtail P \twoheadrightarrow Q$ be an extension. The basic idea of Meyer's proof is to
use the \emph{left} ideal $L \subseteq TP$ that is generated by $S
\subseteq TP$.  Meyer shows
that $L$ is quasi-free and that the complex $X(L)$ is homotopy
equivalent to the kernel of the canonical map $X (TP) \to X (TQ)$.
Secondly, he shows that $X(L)$ is homotopy equivalent to the kernel
of the canonical map $X (TP) \to X(TQ)$. The excision theorem
follows.

\section{Later developments in cyclic theory}

The approach to cyclic theory developed in \cite{CQANs},
\cite{CQCNs}, \cite{CQEx} turned out to be a very natural and
convenient basis for many further advancements in cyclic theory.\mn

The author developed a bivariant $K$-theory for locally convex
algebras and used the excision result of \cite{CQEx} to construct a
bivariant Chern-Connes character from this bivariant theory to
bivariant periodic cyclic theory.\mn

Connes' entire cyclic cohomology \cite{CoEnt} was generalized by
Meyer, using the framework of \cite{CQCNs}, to the bivariant setting
and to bornological algebras, see \cite{Meybook}. Meyer also proved
excision for this theory. \mn

M.Puschnigg developed the very flexible local theory which works for
$C^*$-algebras just as well as for locally convex dense subalgebras
\cite{PuEx}, \cite{PuLoc}.  He obtained the important result that
the local theory of a $C^*$-algebra is independent of the choice of
a natural dense ``smooth'' subalgebra and constructed a bivariant
Chern-Connes character from Kasparov's $KK(A,B)$, for $C^*$-algebras
$A,B$, to his $HC^{loc}(A,B)$.\mn

C.Voigt extended equivariant periodic theory which was known before
only for compact group actions and only for homology \cite{Bry},
\cite{BG}, to actions of general groups and to cyclic cohomology and
to the bivariant theory. A crucial point was again the Hom-complex
for $X$-complexes (which are not quite complexes in this case) for
quasi-free extensions \cite{Voi1}, \cite{Voi2}.

\end{document}